\documentclass[10pt]{article}

\usepackage{a4wide}
\usepackage{amssymb}
\usepackage{amsfonts}
\usepackage{amsmath}
\input xy
\xyoption{arrow} \xyoption{matrix}

\date{}

\newtheorem{proposition}{Proposition}
\newtheorem{theorem}[proposition]{Theorem}

\newtheorem{corollary}[proposition]{Corollary}

\def\GK{{\rm  GK}\,}

\def\der{\partial }

\def\nFM0{{\nu }_{F,M_0}}
\def\nFN0{{\nu }_{F,N_0}}
\def\nGN0{{\nu }_{G,N_0}}

\def\N0{ {\bf N}_0 }

\def\t{\otimes}

\def\ra{\rightarrow}

\def\Xpm{X^{\pm }}

\def\s{\sigma}

\def\l1{{\lambda}_1}

\def\a{\alpha}
\def\a0{ {\alpha }_0}
\def\a1{ {\alpha }_1}

\def\l{\lambda}

\def\nFGM0{{\nu }_{F,G,M_0}}

\def\nFN0{{\nu}_{F,N_0}}

\def\sm{{\sigma}^m}

\def\sm1{{\sigma}^{-1}}

\def\smtp1{{\sigma}^{-t+1}}

\def\S1{S^{-1}}

\def\Xpm1{X^{\pm 1}_1}

\def\sPM1{{\sigma }^{\pm 1}}
\def\sMP1{{\sigma }^{\mp 1 }}

\def\d{\delta}

\def\di{{\rm d.ind}}

\def\L{\Lambda}

\def\CD{{\cal D}}

\def\Ytm1{Y^{t-1}}
\def\Yim1{Y^{i-1}}

\def\Aut{{\rm Aut}}

\def\Der{{\rm Der }}

\def\dim{{\rm dim }}

\def\CJ{ {\cal J}}

\def\SL2Z{ {\rm SL}_2({\bf Z}) }

\def\Gp1{ G^{1 , 1 } }
\def\P11{ P^{-1 , 1 } }
\def\Pp1{ P^{1 , 1 } }

\def\nCLsr{{}^\nu\kern-2pt {\cal L}^{\sigma , \rho  }}
\def\nP{{}^\nu \kern-2pt P}
\def\nL{{}^\nu\kern-2pt L}
\def\nLL{{}^\nu\kern-2pt \Lambda}
\def\nPsr{{}^\nu\kern-2pt P^{\sigma , \rho  }}
\def\nLsr{{}^\nu\kern-2pt L^{\sigma , \rho  }}
\def\nuCL{{}^\nu\kern-2pt  {\cal L}}
\def\nCLsr{{}^\nu\kern-2pt {\cal L}^{\sigma , \rho  }}
\def\nCL1m{{}^\nu\kern-2pt {\cal L}^{-1 , 1  }}
\def\x1nu{x^\frac{1}{\nu}}
\def\xm1nu{x^{-\frac{1}{\nu}}}

\def\ob{\overline{b}}

\def\ra{\rightarrow }

\def\CB{{\cal B}}

\def\nAM0{{\nu }_{{\cal A},M_0}}
\def\nAN0{{\nu }_{{\cal A},N_0}}

\def\End{ {\rm End }}
\def\Der{ {\rm Der }}
\def\CJ{ {\cal J }}

\def\det{ {\rm det }}

\def\di!{\frac{\der^i}{i!}}
\def\dik!{\frac{\der^k_i}{k!}}

\def\hA{\widehat{A}}

\def\hP{\widehat{P}}

\def\N{\mathbb{N}}

\def\0{\overline{0}}
\def\1{\overline{1}}

\def\Ln1{\L_{n,\overline{1}}}

\def\oa{\overline{a}}

\def\a1{a_{\overline{1}}}

\def\S{\Sigma}

\def\vn1{\overrightarrow{n-1}}

\def\hP{\widehat{P}}

\def\mJ{\mathbb{J}}
\def\mI{\mathbb{I}}

\def\K1{{\rm K}_1}

\def\hmI1{\widehat{\mI_1}}
\def\tmI1{\widetilde{\mI_1}}
\def\tmJ1{\widetilde{\mJ_1}}
\def\hB1{\widehat{B_1}}
\def\hCB1{\widehat{\CB_1}}

\def\CB{{\cal B}}

\def\hA{\hat{A}}
\def\hB{\hat{B}}

\begin{document}

\author{V. V.   Bavula 
} 
 
\title{Holonomic modules and  1-generation  in the Jacobian Conjecture}

\maketitle

\begin{abstract} 
A polynomial endomorphism $\s\in \End_K(P_n)$ is called  a {\em Jacobian map} if its Jacobian is a nonzero scalar (the field has zero characteristic). Each Jacobian map $\s$ is extended to an endomorphism $\s$ of the Weyl algebra $A_n$. 

The {\em Jacobian Conjecture} (JC) says that every Jacobian map is an automorphism.  Clearly, the Jacobian Conjecture is true iff the twisted (by $\s$)  $P_n$-module ${}^\s P_n$ is  1-generated for all Jacobian maps $\s$. It is shown that the $A_n$-{\em module}  ${}^\s P_n$ is  1-generated   for all Jacobian maps $\s$. Furthermore, the $A_n$-module  ${}^\s P_n$ is holonomic and as a result has finite length. An explicit upper bound is found for the length of the $A_n$-module  ${}^\s P_n$ in terms of the degree $\deg (\s )$ of the Jacobian map $\s$. Analogous results are given for the Conjecture of Dixmier and the Poisson Conjecture. These results show that the Jacobian Conjecture, the Conjecture of Dixmier and the Poisson Conjecture are questions about holonomic modules for the Weyl algebra $A_n$,  the images of the Jacobian maps,  endomorphisms of the  Weyl algebra $A_n$ and the Poisson endomorphisms are large in the sense that further strengthening of the results on largeness would be either to prove the conjectures or produce  counter examples.

A  short direct algebraic (without reduction to prime characteristic) proof is given of equivalence of the Jacobian and the Poisson Conjectures (this gives a new short proof of equivalence of the Jacobian,   Poisson  and Dixmier Conjectures).  \\




\end{abstract}



In this paper, $K$ is a field of characteristic zero and  $K^\times := K\backslash \{ 0\}$,  $P_n=K[x_1, \ldots , x_n]$ is a polynomial algebra in $n$ the variables, $\Der_K(P_n)$ is the set of all $K$-derivations of the polynomial algebra $P_n$. For a $K$-algebra $A$, the set $\End_K(A)$ is  the monoid of $K$-algebra endomorphisms of $A$ and   $\Aut_K(A)$ is the automorphism group  of $A$. \\

{\bf The  Conjecture of Dixmier, holonomic $A_n$-modules and finite length.} For an endomorphism  $\s\in \End_K (A)$ and an $A$-module  $M$  we denote by  ${}^\s M$  the  $A$-module  $M$  {\bf twisted by $\s$}: ${}^\s M=M$ (as a vector space) and 
$$a\cdot m:= \s (a)m \;\; {\rm for\;all}\;\; a\in A, \;\; m\in M. $$

The ring of differential operators $A_n:=\CD (P_n)$ on the polynomial algebra $P_n$ is called the {\bf Weyl algebra}. The Weyl algebra $A_n$ is generated by the elements $x_1, \ldots , x_n, \der_1, \ldots , \der_n$ subject the defining relations: $[x_i, x_j]=0$, $[\der_i, \der_j]=$ and $[\der_i, x_j]=\d_{ij}$  for all $i,j=1, \ldots , n$ where $\der_i:=\frac{\der}{\der x_i}$, $[a,b]:=ab-ba$, and  $\d_{ij}$ is the Kronecker delta. The Weyl algebra $A_n$ is a simple Noetherian domain of Gelfand-Kirillov dimension $\GK (A_n)=2n.$\\

$\bullet$ The {\bf Inequality of Bernstein}: {\em For all nonzero finitely generated $A_n$-modules} $M$, $$\GK (M)\geq n.$$ 

A finitely generated $A_n$-module $M$ is called {\bf holonomic} if $\GK (M)=n$. Each holonomic module has finite length and is a cyclic $A_n$-module, i.e. 1-generated. 
Each nonzero  sub- or factor module of a holonomic module is holonomic. \\

$\bullet$ The {\bf Conjecture of Dixmier, ${\rm DC}_n$,  \cite{Dix} (1968):} $\End_K(A_n)=\Aut_K(A_n)$.\\

The Weyl algebra $A_n$ is isomorphic to its opposite algebra $A_n^{op}$ via $A_n\ra A_n^{op}$, $x_i\mapsto x_i$, $\der_i\mapsto \der_i$ for $i=1, \ldots , n$. So, the algebra $A_n\t A_n^{op}\simeq A_{2n}$ is isomorphic to the Weyl algebra $A_{2n}$. In particular,  every $A_n$-{\em bimodule} $N$ is a {\em left} $A_{2n}$-module (${}_{A_n}N_{A_n}={}_{A_n\t A_n^{op}}N\simeq {}_{A_{2n}}N$). When we say that an $A_n$-{\em bimodule} $N$ is {\bf holonomic} we mean that the corresponding left $A_{2n}$-module $N$ is holonomic. The Weyl algebra $A_n$ is a simple holonomic $A_n$-bimodule (since $\GK (A_n)=2n$ and $\GK (A_{2n})=4n$).

\begin{theorem}\label{MholsMhol}
\cite[Theorem 1.3]{Bav-RenQn} If  $M$ is a holonomic  $A_n$-module and $\s\in \End_K(A_n)$ then the   $A_n$-module ${}^\s M$ is  also a  holonomic $A_n$-module (and as a result has finite length and is 1-generated over $A_n$).
\end{theorem} 

The Weyl algebra $A_n$ is a simple algebra. So, for each $\s\in \End_K(A_n)$, the image $\s (A_n)$ is isomorphic to the Weyl algebra $A_n$.\\

 $\bullet$  {\em The  Conjecture of Dixmier is true if for every endomorphism $\s\in \End (A_n)$, the $\s(A_n)$-bimodule $A_n$ is  simple.}

\begin{corollary}\label{}
\cite[Corollary 3.4]{Bav-RenQn}  For each algebra endomorphism  $\s :
A_n\ra A_n$, the Weyl algebra $A_n$ is a  holonomic
$\s(A_n)$-bimodule, hence,  of finite length and 1-generated.
\end{corollary}

Each nonzero element $a\in A_n$ is a unique sum $a=\sum_{\alpha, \beta \in \N^n}\l_{\alpha\beta} x^\alpha \der^\beta$ for some scalars $\l_{\alpha\beta}\in K$ where $x^\alpha=x_1^{\alpha_1}\cdots x_n^{\alpha_n}$ and $\der^\beta= \der_1^{\beta_1}\cdots \der_n^{\beta_n}$. The natural number $$\deg (a):=\max \{ |\alpha| +|\beta|\, | \, \l_{\alpha\beta}\neq  0\}$$ is called the {\bf degree} of the element $a$. 
Then $\{ A_{n,i}\}_{i\geq 0}$ is a finite dimensional  filtration of the Weyl algebra $A_n$ where $A_{n,i}:=\{ a\in A_n\,  | \, \deg (a)\leq i\}$ and $\deg (0):=-\infty$ ($A_n=\bigcup_{i\geq 0}A_{n,i}$ and $A_{n,i}A_{n,j}\subseteq A_{n, i+j}$ for all $i,j\geq 0$).

Each  endomorphism $\s \in \End_K(A_n)$ is uniquely determined by the elements 
$$x_1':=\s (x_1), \ldots , x_n':=\s (x_n), \der_i':=\s (\der_1), \ldots , \der_n':=\s (\der_n).$$ 
The natural number
$\deg (\s):=\max\{\deg (x_i'), \deg(\der_i')\, \ \, i=1,\ldots , n \}$ is called the {\bf degree} of $\s$. 

\begin{theorem}\label{A6Nov21}
 Let $\s \in \End(A_n)$  and $d:=\deg (\s)$. Then 
 $L_{\s(A_n)}(A_n)\leq d^{2n} $ where $L_{\s(A_n)}(M)$ is the length of a $\s(A_n)$-bimodule $M$. 
 \end{theorem}
 
  {\it Proof.}  Since $\deg (x_i')\leq d$ and  $\deg (\der_i')\leq d$,
$$ x_i'A_{n,ds}\subseteq A_{n, d(s+1)},\; A_{n,ds} x_i'\subseteq A_{n, d(s+1)},\; \der_i'A_{n,ds}\subseteq A_{n, d(s+1)}\; {\rm and}\;\; A_{n,ds}\der_i'\subseteq A_{n, d(s+1)}\; $$ 
for all $i=1, \ldots , n$ and $s\geq 0$. 
Therefore, $\{ A_{n, ds}\}_{s\geq 0}$ is a finite dimensional filtration of the $\s(A_n)$-bimodule $A_n$ such that 
$$\dim_K(A_{n,ds})={ds+2n\choose 2n}=\frac{1}{(2n)!} (ds+2n)(ds+2n-1)\cdots (ds+1)=\frac{d^{2n}}{(2n)!}s^{2n}+\cdots $$
where three dots denote smaller terms. 
 By \cite[Lemma 8.5.9]{MR},  $L_{\s(A_n)}(A_n)\leq d^{2n} $.  $\Box$\\

{\bf The Jacobian Conjecture, holonomic $A_n$-modules and 1-generation.} Each  endomorphism $\s \in \End_K(P_n)$ is uniquely determined by the polynomials
$$x_1':=\s (x_1), \ldots , x_n':=\s (x_n).$$ 

 The
  matrix of partial derivatives, 
  $$\CJ (\s):=\frac{\der x'}{\der x}:=\bigg(\frac{\der x_i'}{\der x_j}\bigg), \;\; {\rm where}\;\; \CJ (\s)_{ij}:=\frac{\der x_i'}{\der x_j},$$
is called the {\bf Jacobian matrix} of $\s$.  An endomorphism  $\s\in \End (P_n)$ with $\det (\CJ (\s))\in K^\times$ is called a {\bf Jacobian map}. \\

$\bullet$ The {\bf  Jacobian Conjecture, ${\rm JC}_n$  (1939):} {\em Every Jacobian map is an automorphism}.

\begin{theorem}\label{Thm2.1BCW}
\cite[Theorem 2.1]{BCW}  A Jacobian map $\s\in \End (P_n)$ is an automorphism of $P_n$ if the $P_n$-module ${}^\s P_n$ is  finitely generated.
\end{theorem}

$\bullet$ {\em The Jacobian Conjecture is true iff the  $P_n$-module ${}^\s P_n$ is  1-generated for all Jacobian maps} $\s$.\\

 Each Jacobian map $\s$ is  extended to a  (necessarily)  monomorphism of the Weyl algebra $A_n$:  
\begin{equation}\label{sAnext}
\s :A_n\ra A_n , \;\; \der_i\mapsto \der_i', \;\; i=1, \ldots , n,
\end{equation}
where $\der_i'$ is a $K$-derivation of the polynomial algebra $P_n$ which is given by the rule:
\begin{equation}\label{derip}
\der_i'(p):=\frac{1}{\det \, \CJ (\s )}\CJ (\s (x_1), \ldots , \s (x_{i-1}), p, \s(x_{i+1}), \ldots , \s (x_n))\;\; {\rm for\; all}\;\; p\in P_n.
\end{equation}

For an algebra $A$ and its non-empty subset $S$, $C_A(S):=\{ a\in A\, | \, as=sa$ for all $s\in \}$ is the {\em centralizer} of $S$ in $A$. Let $\hP_n:= K[[x_1, \ldots , x_n]]$ and $\hA_n:=\bigoplus_{\alpha\in \N^n}\hP_n\der^\alpha$.  
Proposition \ref{B6Nov21} is a description of all extensions of a Jacobian map of $P_n$  to an endomorphism of the Weyl algebra $A_n$. 

\begin{proposition}\label{B6Nov21}
Let $\s$ be a Jacobian map of $P_n$, $\s$ be its extension to an  endomorphism of the Weyl algebra $A_n$ as in (\ref{sAnext}),  $x_1'=\s(x_1), \ldots , x_n'=\s (x_n)$ and $\der_1'=\s (\der_1), \ldots , \der_n'=\s (\der_n)$, see  (\ref{derip}). 
\begin{enumerate}
\item If $\s'$ is another extension of the Jacobian map $\s$ then $\s'(\der_i)=\der_i'+\der_i'(p)$, $i=1, \ldots , n$ where  $p\in  P_n$, and vice versa.
\item  An extension of the Jacobian map $\s$ of $P_n$ is unique if the images of the elements $\der_1, \ldots , \der_n$ are derivations of $P_n$. So,  the extension $\s$ in (\ref{derip}) is such a unique extension, and $\Der_K(P_n)=\bigoplus_{i=1}^nP_n\der_i'$. 
\item  $C_{A_n}(x_1', \ldots , x_n')=P_n$.
\item Suppose that $x_i'=x_i+\cdots$ for $i=1, \ldots , n$  where the three dots  denote higher terms. Then $\s\in \Aut_K(\hP_n)$ and every extension $\s'$ of the Jacobian map $\s$ to an endomorphism of the Weyl algebra $A_n$ belongs to $\Aut_K(\hA_n)$. 
\end{enumerate}
\end{proposition}

{\it Proof}. 1--3. Up to an affine change of variables in the polynomial algebra $P_n$, we can assume that $\s (x_i)=x_i+\cdots$ for $i=1, \ldots , n$  where the three dots  denote {\em higher} terms. Since $\det (\CJ (\s))\in K^\times$, we have that  $\Der_K(P_n)=\bigoplus_{i=1}^nP_n\der_i'$ (as 
 $\der_i=\sum_{j=1}^n\frac{\der x_j'}{\der x_i}\der_j'$ for all $i=1, \ldots , n$), and so $$\hP_n= K[[x_1', \ldots , x_n']]\;\; {\rm  and}\;\; \s\in \Aut_K(\hP_n).$$ 
If the elements $\s'(\der_i)$ are derivations of the polynomial algebra $P_n$  then $$\s' (\der_i)(x_j')=[\s' (\der_i), x_j']=[\s' (\der_i), \s'(x_j)]=\s' ( [\der_i, x_j])=\s'(\d_{ij})=\d_{ij}\;\; {\rm for}\;\; i,j=1, \ldots , n.$$
Hence,  $\s'(\der_i)=\frac{\der}{\der x_i'}$ for $i=1,\ldots , n$, and so  $\s'(\der_i)=\der_i'$, see (\ref{derip}). 	

In the general case, 
$$ [\s'(\der_i)-\der_i',x_j']=[\s'(\der_i),\s' (x_j)]-[\der_i',x_j']=\d_{ij}-\d_{ij}=0,$$ and so 
$d_i:=\s'(\der_i)-\der_i'\in C_{A_n}(x_1', \ldots , x_n')$. Clearly,
$$ P_n\subseteq C_{A_n}(x_1', \ldots , x_n')\subseteq C_{\hA_n}(x_1', \ldots , x_n')=P_n,$$ and so 
$ C_{A_n}(x_1', \ldots , x_n')=P_n$. Therefore, $d_i\in P_n\subseteq \hP_n= K[[x_1', \ldots , x_n']]$ for all $i=1, \ldots , n$. For all $i,j=1, \ldots, n$, 
 $$0=\s'([\der_i, \der_j])=[\s'(\der_i), \s'(\der_j)]=[\der_i'+d_i, \der_j'+d_j]=\der_i'(d_j)-\der_j'(d_i).$$
 Therefore, there is an element $p\in K[[x_1', \ldots , x_n']]$ such that $d_i=\der_i'(p)$ for $i=1, \ldots , n$, by the Poincar\'{e} Lemma. Since all $d_j\in P_n$, we must have $$\der_i(p)=\sum_{j=1}^n\frac{\der x_j'}{\der x_i}\der_j'(p) =\sum_{j=1}^n\frac{\der x_j'}{\der x_i}d_j\in P_n.$$ Hence, $p\in P_n$ since $K[[x_1',\ldots , x_n']]=\hP_n$.  
 
4. Statement 4 follows from statement 1.    $\Box $

\begin{theorem}\label{6Nov21}
 Let $\s \in \End(P_n)$ be Jacobian map and $d:=\deg (\s)$. Then 
 \begin{enumerate}
\item The $A_n$-module ${}^\s P_n$ is holonomic, hence of finite length and  1-generated as an $A_n$-module.
\item $l_{A_n}({}^{\s}P_n)\leq m^n $ where $m:=\max \{d,  (d-1)^{n-1}-1\}$ where $l_{A_n}(M)$ is the length of an $A_n$-module $M$. 
\end{enumerate}
 \end{theorem}

{\it Proof.} 1. The $A_n$-module $P_n$ is holonomic. By Theorem \ref{MholsMhol},  the $A_n$-module ${}^\s P_n$ is holonomic, hence of finite length and   1-generated as an $A_n$-module. 

2. Since $\deg (x_i')\leq d$, $$\der_i'=\sum_{j=1}^n\frac{\der x_j}{\der x_i'}\der_j=\sum_{j=1}^n\Big(\CJ (\s)^{-1}\Big)_{ij}\der_j, \;\; {\rm   and }\;\; \deg \Big(\CJ (\s)^{-1}\Big)_{ij}\leq  (d-1)^{n-1},$$
we have that  
$$ x_i'P_{n,ms}\subseteq P_{n, m(s+1)},\;\; {\rm and}\;\; \der_i'P_{n,ms}\subseteq P_{n, m(s+1)}\;\; {\rm for\; all}\;\; i=1, \ldots , n\;\; {\rm and}\;\; s\geq 0.$$ 
Therefore, $\{ P_{n, ms}\}_{s\geq 0}$ is a finite dimensional filtration of the $A_n$-module ${}^{\s} P_n$ such that 
$$\dim_K(P_{n,ms})={ms+n\choose n}=\frac{1}{n!} (ms+n)(ms+n-1)\cdots (ms+1)=\frac{m^n}{n!}s^n+\cdots $$
where three dots denote smaller terms. 
By \cite[Lemma 8.5.9]{MR}, $l_{A_n}({}^{\s}P_n)\leq m^n $.  $\Box$\\

The Dixmier Conjecture implies the {\em Jacobian Conjecture}, \cite[page 297]{BCW}),  and the inverse implication
is also true, Tsuchimoto \cite{Tsuchimoto-2005} and Belov-Kanel and Kontsevich \cite{Belov-Konts-2007} (a short proof is given in \cite{Bav-DCJC-2005}).\\

{\bf Equivalence of the Jacobian and  the Poisson Conjectures.} The Weyl algebra $A_n=\CD (P_n)=\bigcup_{i\geq 0}\CD (P_n)_i$ is a ring of differential operators on $P_n$ and hence admits the {\bf degree filtration} $\{ \CD (P_n)_i\}_{i\geq 0}$ where $\CD (P_n)_i=\bigoplus_{\{ \alpha\in \N^n\, | \, |\alpha |\leq i\}}P_n\der^\alpha$. The {\bf associated graded algebra} $${\rm gr}(A_n):=\bigoplus_{i\geq 0}{\rm gr}(A_n)_i,$$  where ${\rm gr}(A_n)_i:=\CD (P_n)_i/\CD (P_n)_{i-1}$ and $\CD (P_n)_{-1}:=0$,  is a polynomial algebra $P_{2n}$ in $2n $ variables $x_1, \ldots , x_n,  x_{n+1}, \ldots ,  x_{n+n}$ (where $ x_{n+i}:=\der_i+P_n$) that admits the canonical Poisson structure given by the rule: 
\begin{equation}\label{grAnPois}
\{\cdot, \cdot \}: {\rm gr}(A_n)_i\t_K {\rm gr}(A_n)_j\ra {\rm gr}(A_n)_{i+j-1}, \;\; (\oa, \ob)\mapsto \{ \oa, \ob\}:=[a,b]+\CD (P_n)_{i+j-2}
\end{equation}
  where $\oa :=a+\CD (P_n)_{i-1}$ and $\ob :=b+\CD (P_n)_{j-1}$ (since $[\CD (P_n)_i,\CD (P_n)_j]\subseteq\CD (P_n)_{i+j-1}$ for all $i,j\geq 0$). Equivalently,
\begin{equation}\label{grAnPois1}
\{x_i, x_j\}=0,\;\;  \{ x_{n+i},  x_{n+j}\}=0\;\; {\rm and}\;\; \{  x_{n+i}, x_j\}= \d_{ij}\;\; {\rm for\;all}\;\; i,j=1, \ldots ,n .
\end{equation}

$\bullet$ {\bf The Poisson Conjecture, ${\rm PC}_{2n}$:} $\End_{\rm Pois}(P_{2n})=\Aut_{\rm Pois}(P_{2n})$. \\

The Poisson Conjecture and the Conjecture of Dixmier are equivalent (Adjamagbo and  van den Essen, \cite{Adjam-vdEs-2007}).

\begin{theorem}\label{}\marginpar{}
\begin{enumerate}
\item ${\rm JC}_{2n}\Rightarrow {\rm PC}_{2n}$.
\item ${\rm DC}_{2n}\Rightarrow {\rm PC}_{2n}$.
\item ${\rm PC}_{2n}\Rightarrow {\rm JC}_n$.  
\item The Jacobian Conjecture and the Poisson Conjecture are equivalent. 
\item The Jacobian Conjecture, the Conjecture of Dixmier and Poisson Conjecture are equivalent. 
\end{enumerate}
\end{theorem}

{\it Proof}. 1. Given  $\s\in \End_{\rm Pois}(P_{2n})$. Then $\det(\CJ (\s))\in \{ \pm 1\}$ (see the proof of Step 6 of \cite[Theorem 3]{Bav-DCJC-2005} of the fact that ${\rm JC}_{2n}\Rightarrow {\rm DC}_n$): Notice that $\det(\{ x_i, x_j\})\in \{ \pm 1\}$ where $1\leq i,j\leq 2n$, and so 
\begin{eqnarray*}
\{ \pm 1\}&\ni &\det(\{ x_i, x_j\})=\s(\det(\{ x_i, x_j \} ))=  \det(\s (\{ x_i, x_j\}))\\
&=&  \det( \{ \s(x_i), \s (x_j)\})=\det (\CJ^t(\s)\cdot \Big(\{ x_i, x_j\}\Big)\cdot \CJ (\s))\\
& =&\det(\CJ(\s))^2\det(\{ x_i, x_j\}),
\end{eqnarray*}
and so $\det(\CJ (\s))\in \{ \pm 1\}$. By ${\rm JC}_{2n}$, $\s\in \Aut_K(P_{2n})$,  and statement 1  follows. 

2. Given  $\s\in \End_{\rm Pois}(P_{2n})$. Then 
maps 
\begin{equation}\label{sPoisWA}
\s : A_{2n} \ra A_{2n}, \;\; x_i\mapsto x_i':=
\s (x_i), \;\; \der_i\mapsto \der_i'(\cdot ):=\begin{cases}
\{x_{n+i}', \cdot \} & \text{if }i=1,\ldots , n,\\
\{-x_{n-i}', \cdot \}& \text{if }i=n+1,\ldots , 2n,\\
\end{cases}
\end{equation}
is an algebra endomorphism of the Weyl algebra $A_{2n}$ where $\der_i'\in \Der_K(P_{2n})$. By  ${\rm DC}_{2n}$, $\s\in \Aut_K(A_{2n})$, and so $$A_{2n}=\bigoplus_{\alpha, \beta \in \N^{2n}}Kx'^\alpha \der'^\beta.$$ It follows that $\CD (P_{2n})_i=\bigoplus_{ \{ \alpha , \beta \in \N^{2n}\, | \, |\beta|\leq i\}}Kx'^\alpha \der'^\beta $ (use the defining relations in the new variables of $A_{2n}$). By the very definition,  the automorphism $\s$ respects the degree filtration on $A_{2n}$. In particular, $\s (P_{2n})=P_{2n}$ since $\CD (P_{2n})_0=P_{2n}$, i.e. $\s \in \Aut_{\rm Pois}(P_{2n})$. 

3. Given a Jacobian map $\s \in \End_K(P_n)$. Let $\s\in \End_K(A_n)$  be its  extension given by (\ref{sAnext}). 
By (\ref{derip}), 
$$\s (\CD (P_n)_i)\subseteq \CD (P_n)_i\;\; {\rm for\; all}\;\; i\geq 0,$$ i.e. the endomorphism $\s$ of the Weyl algebra $A_n$ respects the degree filtration and so the associated graded map 
\begin{equation}\label{grsPois}
{\rm gr}(\s): {\rm gr}(A_n)\ra {\rm gr}(A_n), \;\; a+\CD (P_n)_{i-1}\mapsto \s (a)+\CD (P_n)_{i-1}
\end{equation}
respects the Poisson structure, i.e. ${\rm gr}(\s) \in \End_{\rm Pois}({\rm gr}(A_n))$. By ${\rm PC}_{2n}$, ${\rm gr}(\s) \in \Aut_{\rm Pois}({\rm gr}(A_n))$, hence $\s\in \Aut_K(P_n)$ since the automorphism ${\rm gr}(\s) \in \Aut_{\rm Pois}({\rm gr}(A_n))$ is a {\em graded} automorphism and $\CD (P_n)_0=P_n$.

4. Statement  4 follows from statements 1 and 3. 

5. Statement 5 follows from statement 4 and the equivalence ${\rm JC}_{2n}\Leftrightarrow {\rm DC}_n$. $\Box $

By (\ref{sPoisWA}),\\

$\bullet$ ${\rm PC}_{2n}$ {\em is true iff the $A_{2n}$-module ${}^{\s}P_{2n}$ is simple for all $\s\in \End_{\rm Pois}(P_{2n})$.}

\begin{theorem}\label{16Nov21}
 Let $\s \in \End_{\rm Pois}(P_{2n})$ and $d:=\deg (\s)$. Then 
 \begin{enumerate}
\item The $A_{2n}$-module ${}^\s P_{2n}$ is holonomic, hence of finite length and  1-generated as an $A_{2n}$-module.
\item $l_{A_{2n}}({}^{\s}P_{2n})\leq d^{2n} $. 
\end{enumerate}
 \end{theorem}

{\it Proof.} 1. The $A_{2n}$-module $P_{2n}$ is holonomic. By Theorem \ref{MholsMhol},  the $A_{2n}$-module ${}^\s P_{2n}$ is holonomic, hence of finite length and   1-generated as an $A_{2n}$-module. 

2. Since $\deg (x_i')\leq d$ and $\deg (\der_i')\leq d$ (see (\ref{sPoisWA})),
$$ x_i'P_{2n,ds}\subseteq P_{2n, d(s+1)},\;\; {\rm and}\;\; \der_i'P_{2n,ds}\subseteq P_{2n, d(s+1)}\;\; {\rm for\; all}\;\; i=1, \ldots , 2n\;\; {\rm and}\;\; s\geq 0.$$ 
Therefore, $\{ P_{2n, ds}\}_{s\geq 0}$ is a finite dimensional filtration of the $A_{2n}$-module ${}^{\s} P_{2n}$ such that 
$$\dim_K(P_{2n,ds})={ds+2n\choose 2n}=\frac{1}{(2n)!} (ds+2n)(ds+2n-1)\cdots (ds+1)=\frac{d^{2n}}{(2n)!}s^{2n}+\cdots $$
where three dots denote smaller terms. 
By \cite[Lemma 8.5.9]{MR}, $l_{A_{2n}}({}^{\s}P_{2n})\leq d^{2n} $.  $\Box$\\

{\bf An analogue of the Conjecture of Dixmier for the algebras $\mI_n$ of integro-differential operators.} Let  
$\mI_n:=K\langle x_1, \ldots , x_n, \der_1, \ldots
,\der_n, \int_1, \ldots , \int_n\rangle $ be the   algebra of
polynomial  integro-differential operators where $\int_i: P_n\ra P_n$, $ p\mapsto \int p \, dx_i$, i.e. $\int_i x^\alpha = (\alpha_i+1)^{-1}x_ix^\alpha$ for all $\alpha \in \N^n$, \cite{intdifaut}.

$\bullet$ {\bf Conjecture, \cite{IntDifDixConj} (2012):} $\End_K(\mI_n)=\Aut_K(\mI_n)$. \\

\begin{theorem}\label{}\marginpar{}
 \cite[Theorem 1.1]{IntDifDixConj}  $\End_K(\mI_1)=\Aut_K(\mI_1)$.
\end{theorem}

{\bf An analogue of the Jacobian Conjecture and the Conjecture of Dixmier for the algebras $A_{n,m}:=A_n\t P_m$.}
 The centre of the algebra $A_{n,m}$ is  $P_m$. Hence, for all $\s\in \Aut_K(A_{n,m})$, $\s (P_m)=P_m$.\\

$\bullet$  {\bf  Conjecture, \cite{Bav-InvForWeylAlgPol-2007} (2007),  ${\rm JD}_{n,m}$}:  
{\em Every endomorphism $\s :A_n\t P_m\ra
A_n\t P_m$ such that  $\s (P_m)\subseteq P_m$ and
$\det \bigg(\frac{\der \s (x_i)}{\der x_j}\bigg)\in K^*$ is  an  automorphism.}\\

\begin{theorem}\label{}\marginpar{}
\cite[Theorem 5.8, Proposition 5.9]{Bav-InvForWeylAlgPol-2007} ${\rm JD}_{n,m} \Leftrightarrow {\rm JC}_m +{\rm DC}_n$.
\end{theorem}

$${\bf Acnowledgements} $$

The author would like to thank H. Bass and D. Wright for the recent discussions on the Jacobian Conjecture (November 2021) and the Royal Society for support.

School of Mathematics and Statistics

University of Sheffield

Hicks Building

Sheffield S3 7RH

UK

email: v.bavula@sheffield.ac.uk

\end{document}